\documentclass{article}

\usepackage{amssymb}
\usepackage{amsmath}

\usepackage[dvips]{graphicx}

\begin{document}

 
\begin{center}
\textbf{ Functional continuity of unital $B_{0}$-algebras with orthogonal bases}                                                                                                                                                                                   
\end{center}

\noindent \textbf{}

\begin{center}
M. El Azhari
\end{center}

\noindent \textbf{ } 

\noindent \textbf{Abstract.} Let $A$ be a unital $B_{0}$-algebra with an orthogonal basis, then every multiplicative linear functional on $A$ is continuous. This gives an answer to a problem posed by Z. Sawon and Z. Wronski.   

\noindent \textbf{}



\noindent \textbf{}
 
\noindent \textbf{} 
 
\noindent\textbf{I. Preliminaries}

\noindent \textbf{} 

\noindent \textbf{} A topological algebra is a complex associative algebra which is also a Hausdorff topological vector space such that the multiplication is separately continuous. A locally convex algebra is a topological algebra whose topology is determined by a family of seminorms. A complete metrizable locally convex algebra is called a $B_{0}$-algebra. The topology of a $B_{0}$-algebra $A$ may be given by a countable family $(\Vert .\Vert_{i})_{i\geq 1}$ of seminorms such that $\Vert x\Vert_{i}\leq\Vert x\Vert_{i+1}$ and $\Vert xy\Vert_{i}\leq\Vert x\Vert_{i+1}\Vert y\Vert_{i+1}$ for all $i\geq 1$ and $x, y \in A.$ A multiplicative linear functional on a complex algebra $A$ is an algebra homomorphism from $A$ to the complex field. Let $A$ be a topological algebra. $M^{\ast}(A)$ denotes the set of all nonzero multiplicative linear functionals on $A$. $M(A)$ denotes the set of all nonzero continuous multiplicative linear functionals on $A.$ A seminorm $p$ on $A$ is lower semicontinuous if the set $\lbrace x\in A: p(x)\leq 1\rbrace$ is closed in  $A.$

\noindent \textbf{} Let $A$ be a topological algebra. A sequence $(e_{n})_{n\geq 1}$ in $A$ is a basis if for each $x\in A$ there is a unique sequence $(\alpha_{n})_{n\geq 1}$ of complex numbers such that $x = \Sigma_{n=1}^{\infty}\alpha_{n}e_{n}.$ Each linear functional $e_{n}^{\ast}:A\rightarrow \mathbb{C},\, e_{n}^{\ast}(x)=\alpha_{n},$ is called a coefficient functional. If each $e_{n}^{\ast}$ is continuous, the basis $(e_{n})_{n\geq 1}$ is called a Schauder basis. A basis $(e_{n})_{n\geq 1}$ is orthogonal if $e_{i}e_{j}=\delta_{ij}e_{i}$ where $\delta_{ij}$ is the Kronecker symbol. If $(e_{n})_{n\geq 1}$ is an orthogonal basis, then each $e_{n}^{\ast}$ is a multiplicative linear functional on $A.$ Let $A$ be a topological algebra with an orthogonal basis $(e_{n})_{n\geq 1}.$ If $A$ has a unity $e,$ then $e = \Sigma_{n=1}^{\infty}e_{n}.$ Let $(x_{k})_{k}$ be a net in $A$
converging to $0,$ since the multiplication is separately continuous, $e_{n}x_{k}= e_{n}^{\ast}(x_{k})e_{n}\rightarrow_{k} 0$ and so $e_{n}^{\ast}(x_{k})\rightarrow_{k} 0.$ Then each orthogonal basis in a topological algebra is a Schauder basis. Let $f$ be a multiplicative linear functional on $A.$ If 
$f(e_{n_{0}})\neq 0$ for some $n_{0}\geq 1,$ then $f(x)f(e_{n_{0}})= f(xe_{n_{0}})= f(e_{n_{0}}^{\ast}(x)e_{n_{0}})= e_{n_{0}}^{\ast}(x)f(e_{n_{0}})$ for all $x\in A$ and therefore $f = e_{n_{0}}^{\ast}\in M(A).$ This shows that
$M(A)=\lbrace e_{n}^{\ast}: n\geq 1\rbrace.$

\noindent \textbf{} Here we consider unital $B_{0}$-algebras with orthogonal bases. These algebras were investigated in [3] where we can find examples of such algebras.

\noindent \textbf{} 

\noindent \textbf{II. Results} 

\noindent \textbf{}

\noindent \textbf{Proposition II.1.} Let $(A,(\Vert .\Vert_{i})_{i\geq 1})$ be a unital $B_{0}$-algebra with an orthogonal basis $(e_{n})_{n\geq 1}.$ Then there exists $x\in A$ such that the sequence $(e_{n}^{\ast}(x))_{n\geq 1}$ is not bounded.

\noindent \textbf{}  

\noindent \textbf{Proof.} Suppose that $ \sup _{n\geq 1}\vert e_{n}^{\ast}(x)\vert <\infty$ for all     $x\in A.$ For each $x\in A,$ let $\Vert x\Vert = \sup _{n\geq 1}\vert e_{n}^{\ast}(x)\vert$, $\Vert .\Vert  $ is a lower semicontinuous norm on $A$ since $(e_{n})_{n\geq 1}$ is a Schauder basis. Let $\tau_{A}$ be the topology on $A$ determined by the family $(\Vert .\Vert_{i})_{i\geq 1}$ of seminorms. We define a new topology $\tau$ on $A$ described by the norm $\Vert .\Vert$ and the family $(\Vert .\Vert_{i})_{i\geq 1}$ of seminorms. The topology $\tau$ is stronger than the topology $\tau_{A}.$ By Garling's completeness theorem [1, Theorem 1], $(A,\tau)$ is complete. The topologies $\tau_{A}$ and $\tau$ are homeomorphic by the open mapping theorem. Then there exist $i_{0}\geq 1$ and $M > 0$ such that $\Vert x\Vert\leq M\Vert x\Vert_{i_{0}}$ for all $x\in A,$ hence $1 = \vert e_{n}^{\ast}(e_{n})\vert\leq\Vert e_{n}\Vert\leq M\Vert e_{n}\Vert_{i_{0}}$ for all $n\geq 1.$ This contradicts the fact that $e_{n}\rightarrow_{n} 0.$

\noindent \textbf{}

\noindent \textbf{Proposition II.2.} Let $A$ be a unital $B_{0}$-algebra with an orthogonal basis $(e_{n})_{n\geq 1}.$ If $x = \Sigma_{n=1}^{\infty}t_{n}e_{n}\in A  $ such that $t_{n}\in \mathbb{R},\,t_{n}\leq t_{n+1}$ for $n\geq 1$ and $t_{n}\rightarrow\infty,$ then $f(x)\in \mathbb{R}\,$ for all $f\in M^{\ast}(A).$

\noindent \textbf{}

\noindent \textbf{Proof.} If $f\in M(A),\, f = e_{n}^{\ast}$ for some $n\geq 1,$ then $f(x) = t_{n}\in \mathbb{R}.  $ If $f\in M^{\ast}(A)\smallsetminus M(A),$ then $f(e_{n}) = 0$ for all $n\geq 1.$ Suppose that $f(x)\notin \mathbb{R},\, f(x) = \alpha + i\beta$ with $\beta \neq 0.$ Since $t_{n}\rightarrow\infty,$ there exists $n_{0}\geq 1  $ such that $ t_{n}\geq\alpha +\vert\beta\vert$ for all $n\geq n_{0}.$ We define the sequence $(s_{n})_{n\geq 1}$ by $s_{n} = t_{n_{0}}$ for $1\leq n\leq n_{0}$ and $s_{n} = t_{n}$ for $n\geq n_{0}+1.  $ It is clear that $y = \Sigma_{n=1}^{\infty}s_{n}e_{n}\in A$ such that $s_{n}\geq\alpha +\vert\beta\vert,\, s_{n}\leq s_{n+1}$ for $n\geq 1$ and $s_{n}\rightarrow \infty.$ Since $f(e_{n})= 0$ for all $n\geq 1,\,f(y)= f(x)= \alpha + i\beta.$ We have $f(\vert\beta\vert^{-1}y)= \vert\beta\vert^{-1}\alpha + i\vert\beta\vert^{-1}\beta,$ then $f(\vert\beta\vert^{-1}y - \vert\beta\vert^{-1}\alpha e)= i\vert\beta\vert^{-1}\beta.$ Set $z = \vert\beta\vert^{-1}y - \vert\beta\vert^{-1}\alpha e = \Sigma_{n=1}^{\infty}\frac{s_{n}-\alpha}{\vert\beta\vert}e_{n}.$ The real sequence $(\frac{s_{n}-\alpha}{\vert\beta\vert})_{n\geq 1}$ is positive increasing and $\frac{s_{n}-\alpha}{\vert\beta\vert}\rightarrow \infty,$ then $z^{-1}= \Sigma_{n=1}^{\infty}\frac{\vert\beta\vert}{s_{n}-\alpha}e_{n}\in A$ by [3, Theorem 0.1] and $f(z^{-1})= -i\vert\beta\vert\beta^{-1},$ so $f(z + z^{-1})= 0.$ Set $v = z + z^{-1}= \Sigma_{n=1}^{\infty}(\frac{s_{n}-\alpha}{\vert\beta\vert}+ \frac{\vert\beta\vert}{s_{n}-\alpha})e_{n}   $ and $v_{n}= \frac{s_{n}-\alpha}{\vert\beta\vert}+ \frac{\vert\beta\vert}{s_{n}-\alpha}$ for all $n\geq 1.$ Since the map $g: [1,\infty)\rightarrow \mathbb{R},\,g(x) = x + \frac{1}{x},$ is increasing and the sequence 
$(\frac{s_{n}-\alpha}{\vert\beta\vert})_{n\geq 1}\subset [1,\infty)$ is increasing, it follows that $(v_{n})_{n\geq 1}$ is a positive increasing sequence and $v_{n}\rightarrow\infty.$ By [3, Theorem 0.1], 
$v^{-1}= \Sigma_{n=1}^{\infty}\frac{1}{v_{n}}e_{n}\in A$ and therefore $f(e) = f(v)f(v^{-1})= 0.$ This contradicts the fact that $f$ is nonzero.

\noindent \textbf{}

\noindent \textbf{} The following two results are due to Sawon and Wronski [3], the proofs are given for completeness.

\noindent \textbf{}

\noindent \textbf{Theorem II.3}([3, Theorem 2.1]). Let $A$ be a unital $B_{0}$-algebra with an orthogonal basis $(e_{n})_{n\geq 1}.$ If $x = \Sigma_{n=1}^{\infty}t_{n}e_{n}\in A  $ such that $t_{n}\in \mathbb{R},\,t_{n}\leq t_{n+1}$ for $n\geq 1$ and $t_{n}\rightarrow\infty,$ then every multiplicative linear functional on $A$ is continuous. 

\noindent \textbf{}

\noindent \textbf{Proof.} Suppose that $M^{\ast}(A)\smallsetminus M(A)$ is nonempty. Let $f\in M^{\ast}(A)\smallsetminus M(A),$ then $f(e_{n})= 0\,$ for all $n\geq 1.$ Put $f(x)=\alpha,\,$ then $\alpha\in\mathbb{R}$ by Proposition II.2. Since $t_{n}\rightarrow\infty,$ there exists $n_{0}\geq 1$ such that $t_{n} > \alpha$ for $n\geq n_{0}.$ Consider  $y = \Sigma_{n=1}^{\infty}\lambda_{n}e_{n}\in A$ such that $\lambda_{n}= t_{n_{0}}$ for $1\leq n\leq n_{0}$ and $\lambda_{n}= t_{n}$ for $n\geq n_{0}+1.$ Since $f(e_{n})= 0 $ for all $n\geq 1,$ it follows that $f(x)= f(y)=\alpha.$ We have                                               $y -\alpha e = \Sigma_{n=1}^{\infty}v_{n}e_{n}\in A$ where $v_{n}=\lambda_{n} -\alpha\,$ for all $n\geq 1.$ It is clear that $v_{n} > 0,\,v_{n}\leq v_{n+1}$ for $n\geq 1$ and $v_{n}\rightarrow\infty.$ By [3, Theorem 0.1],     $(y -\alpha e)^{-1}= \Sigma_{n=1}^{\infty}\frac{1}{v_{n}}e_{n}\in A,$ so $y -\alpha e$ is invertible and        $f(y -\alpha e)= 0,$ a contradiction.

\noindent \textbf{}

\noindent \textbf{Proposition II.4}([3, p.109]). Let $(A,(\Vert .\Vert_{i})_{i\geq 1})$ be a unital $B_{0}$-algebra with an orthogonal basis $(e_{n})_{n\geq 1}.$ Then the set $N$ of all positive integers can be split into two disjoint subsets $N_{1}$ and $N_{2}$ such that by putting $A_{1}= \overline{span}(e_{n})_{n\in N_{1}}$ and $A_{2}= \overline{span}(e_{n})_{n\in N_{2}},$ we have\\
(1) $A = A_{1}\oplus A_{2};$\\
(2) if $f$ is a multiplicative linear functional on $A$ such that $f\notin M(A),$ then $f_{/A_{1}}= 0.$  
 
\noindent \textbf{} 

\noindent \textbf{Proof.} By Proposition II.1, there is $x = \Sigma_{n=1}^{\infty}t_{n}e_{n}\in A$ such that the sequence $(t_{n})_{n\geq 1}$ is not bounded. Then there exists a subsequence $(t_{k_{n}})_{n\geq 1}$ of $(t_{n})_{n\geq 1}$ such that $\vert t_{k_{n}}\vert\geq n^{2}$ for all $n\geq 1.$ For each $i\geq 1,$ there is $M_{i} > 0$ such that $\Vert t_{n}e_{n}\Vert_{i}\leq M_{i}$ for all $n\geq 1.$ Let $i\geq 1$ and $n\geq 1,\, n^{2}\Vert e_{k_{n}}\Vert_{i}\leq\vert t_{k_{n}}\vert\Vert e_{k_{n}}\Vert_{i}= \Vert t_{k_{n}}e_{k_{n}}\Vert_{i}\leq M_{i},$ then $\Vert e_{k_{n}}\Vert_{i}\leq n^{-2}M_{i}.$ This implies that
$\Sigma_{n=1}^{\infty}e_{k_{n}}$ is absolutely convergent. Let $A_{1}= \overline{span}\lbrace e_{k_{n}}: n\geq 1\rbrace,\,A_{1}$ is a unital $B_{0}$-algebra with an orthogonal basis $(e_{k_{n}})_{n\geq 1}$ and $\Sigma_{n=1}^{\infty}n^{\frac{1}{2}}e_{k_{n}}\in A_{1}$ since $n^{\frac{1}{2}}\Vert e_{k_{n}}\Vert_{i}\leq n^{-\frac{3}{2}}M_{i}$ for all $i\geq 1$ and $n\geq 1.$ Set $N_{1}=\lbrace k_{n}: n\geq 1\rbrace,\,N_{2}= N\smallsetminus N_{1}  $ and $A_{2}=\overline{span}\lbrace e_{n}: n\in N_{2}\rbrace.\,A_{2}$ is a $B_{0}$-algebra with an orthogonal basis $(e_{n})_{n\in N_{2}}$ and the unity $u_{2}= e - u_{1}$ where $u_{1}=\Sigma_{n=1}^{\infty}e_{k_{n}}$ is the unity of $A_{1}.$ Let $x\in A,\,x = xe = x(u_{1}+u_{2}) = xu_{1}+ xu_{2}\in A_{1}+ A_{2},$ then $A = A_{1}\oplus A_{2}.$ If $f$ is a multiplicative linear functional on $A$ such that $f\notin M(A),\,f_{/A_{1}}$ is a multiplicative linear functional on $A_{1}$ such that $f_{/A_{1}}(e_{k_{n}})= 0$ for all $n\geq 1.$ Since $ \Sigma_{n=1}^{\infty}n^{\frac{1}{2}}e_{k_{n}}\in A_{1},\,f_{/A_{1}}$ is continuous on $A_{1}$ by    Theorem II.3 and therefore $f_{/A_{1}}= 0.$

\noindent \textbf{}

\noindent \textbf{} Sawon and Wronski [3, p.109] posed the following problem:

\noindent \textbf{}

\noindent \textbf{Problem.}  Let $A$ be a unital $B_{0}$-algebra with an orthogonal basis $(e_{n})_{n\geq 1}.$ Does there exist a maximal subalgebra $A_{1}^{'}=\overline{span}\lbrace e_{n}: n\in N_{1}^{'}\rbrace (N_{1}^{'}\subset N)$ of $A$ for which (1) and (2) hold?

\noindent \textbf{}

\noindent \textbf{Proposition II.5.} Let $A$ be a unital $B_{0}$-algebra with an orthogonal basis $(e_{n})_{n\geq 1}.$ Then the following assertions are equivalent:\\
$(i)\,A = \overline{span}\lbrace e_{n}: n\in N\rbrace$ is a maximal subalgebra of itself for which (1) and (2) hold;\\
$(ii)$ every multiplicative linear functional on $A$ is continuous.

\noindent \textbf{}

\noindent \textbf{Proof.} $(i)\Rightarrow (ii):$ Let $f$ be a multiplicative linear functional on $A$ such that $f\notin M(A),\,f$ is zero on $A$ by $(i).$ Then every multiplicative linear functional on $A$ is continuous.\\
$(ii)\Rightarrow (i):$ It is clear that $A$ satisfies (1). Let $f$ be a multiplicative linear functional on $A$ such that $f\notin M(A),$ then $f$ is zero on $A$ since $f$ is continuous, hence $A$ satisfies (2).

\noindent \textbf{}

\noindent \textbf{Proposition II.6.} Let $(t_{n})_{n\geq n_{0}}$ be a complex sequence, the following assertions are equivalent:\\
$(i)\,\Sigma_{n=n_{0}}^{\infty}\vert t_{n}- t_{n+1}\vert < \infty; $\\
$(ii)$ there exists $M > 0$ such that $\vert t_{q}\vert + \Sigma_{n=p}^{q-1}\vert t_{n}- t_{n+1}\vert\leq M$ for all $q > p \geq n_{0}.$ 

\noindent \textbf{}

\noindent \textbf{Proof.} $(i)\Rightarrow (ii):$ Let $\epsilon > 0,$ there exists $n_{1}\geq n_{0}$ such that $\Sigma_{k=n}^{\infty}\vert t_{k}- t_{k+1}\vert\leq\epsilon$ for every $n\geq n_{1}.$ let $m > n\geq n_{1},\,\vert t_{n}- t_{m}\vert\leq\vert t_{n}- t_{n+1}\vert +...+\vert t_{m-1}- t_{m}\vert\leq\epsilon.$
Then the sequence $(t_{n})_{n\geq n_{0}}$ converges, so there is $M_{0} > 0$ such that $\vert t_{n}\vert\leq M_{0}$ for all $n\geq n_{0}.$ Let $q > p\geq n_{0},\,\vert t_{q}\vert + \Sigma_{n=p}^{q-1}\vert t_{n}- t_{n+1}\vert\leq M_{0}+ \Sigma_{n=n_{0}}^{\infty}\vert t_{n}- t_{n+1}\vert.
$\\
$(ii)\Rightarrow (i):$ Let $p\geq n_{0}$ and $q = p + 1,\,\Sigma_{n=n_{0}}^{p}\vert t_{n}- t_{n+1}\vert\leq\vert t_{q}\vert + \Sigma_{n=n_{0}}^{q-1}\vert t_{n}- t_{n+1}\vert\leq M.$ Then the sequence $(\Sigma_{n=n_{0}}^{p}\vert t_{n}- t_{n+1}\vert)_{p\geq n_{0}}$ is positive increasing and bounded, so it is convergent i.e. $\Sigma_{n=n_{0}}^{\infty}\vert t_{n}- t_{n+1}\vert < \infty. $

\noindent \textbf{}

\noindent \textbf{Proposition II.7.} Let $(A,(\Vert .\Vert_{i})_{i\geq 1})$ be a unital $B_{0}$-algebra with an orthogonal basis $(e_{n})_{n\geq 1}.$ If $(t_{n})_{n\geq n_{0}}$ is a complex sequence such that
$\Sigma_{n=n_{0}}^{\infty}\vert t_{n}- t_{n+1}\vert < \infty,$ then $\Sigma_{n=n_{0}}^{\infty}t_{n}e_{n}\in A.$

\noindent \textbf{}

\noindent \textbf{Proof.} Let $q > p\geq n_{0},$ by using the equality $t_{n}= t_{q}+ \Sigma_{k=n}^{q-1}(t_{k}- t_{k+1})$ for every $p\leq n < q,$ we obtain that $\Sigma_{n=p}^{q-1}t_{n}e_{n}= t_{q}(e_{p}+...+e_{q-1})+ \Sigma_{k=p}^{q-1}(t_{k}- t_{k+1})(e_{p}+...+e_{k}).$ Let $i\geq 1,\,\Vert\Sigma_{n=p}^{q-1}t_{n}e_{n}\Vert_{i}\leq \vert t_{q}\vert\Vert e_{p}+...+e_{q-1}\Vert_{i}+\Sigma_{k=p}^{q-1}\vert t_{k}- t_{k+1}\vert\Vert e_{p}+...+e_{k}\Vert_{i}\leq (\vert t_{q}\vert+ \Sigma_{k=p}^{q-1}\vert t_{k}- t_{k+1}\vert)\sup _{p\leq k\leq q}\Vert e_{p}+...+e_{k}\Vert_{i}\leq M\sup _{p\leq k\leq q}\Vert e_{p}+...+e_{k}\Vert_{i}$ by Proposition II.6. Let $\epsilon > 0,$ since  
$e = \Sigma_{n=1}^{\infty}e_{n}\in A,$ there is $n_{1}\geq n_{0}$ such that $\Vert e_{p}+...+e_{k}\Vert_{i}\leq \epsilon M^{-1}$ for $n_{1}\leq p\leq k,$ hence $\,\sup _{p\leq k\leq q}\Vert e_{p}+...+e_{k}\Vert_{i}\leq \epsilon M^{-1}$ for $n_{1}\leq p < q.$ Then $\Vert \Sigma_{n=p}^{q-1}t_{n}e_{n}\Vert_{i}\leq \epsilon$ for $n_{1}\leq p < q.$ This shows that $\Sigma_{n=n_{0}}^{\infty}t_{n}e_{n}$ is convergent in $A.$

\noindent \textbf{}

\noindent \textbf{Theorem II.8.} Let $(A,(\Vert .\Vert_{i})_{i\geq 1})$ be a unital $B_{0}$-algebra with an orthogonal basis $(e_{n})_{n\geq 1}.$ Then every multiplicative linear functional on $A$ is continuous.

\noindent \textbf{}

\noindent \textbf{Proof.} By Proposition II.7, $x = \Sigma_{n=1}^{\infty}\frac{1}{n}e_{n}\in A.$ Let $f\in M^{\ast}(A)$ and $\alpha = f(x),$ then $f(\alpha e - x)= 0.$ We have $\alpha e - x = \Sigma_{n=1}^{\infty}(\alpha - \frac{1}{n})e_{n}= \Sigma_{n=1}^{\infty}\frac{\alpha n - 1}{n}e_{n}.$ Let $\mathbb{N}$ be the set of all positive integers. Put $I_{\alpha}=\lbrace n\in \mathbb{N}: n\neq \alpha^{-1}\rbrace,\,I_{\alpha}= \mathbb{N}$ if $\alpha^{-1}\notin \mathbb{N}$ and $I_{\alpha}= \mathbb{N}\smallsetminus\lbrace\alpha^{-1}\rbrace$ if $\alpha^{-1}\in \mathbb{N}.$ Let $m_{\alpha}=\inf\lbrace n\in \mathbb{N}: \vert \alpha\vert n - 1 > 0\rbrace$ and consider the complex sequence $(\frac{n}{\alpha n - 1})_{n\geq m_{\alpha}}.$ Let $n\geq m_{\alpha},\,\,\frac{n}{\alpha n - 1}-\frac{n+1}{\alpha (n+1)- 1}=\frac{1}{(\alpha n - 1)(\alpha(n+1) - 1)}.$ We have $\vert\alpha n - 1\vert\geq\vert\alpha\vert n -1$ and $\vert\alpha (n + 1) - 1\vert\geq\vert\alpha\vert(n + 1)- 1 = \vert\alpha\vert n + \vert\alpha\vert - 1\geq\vert\alpha\vert n - 1.$ Let $n\geq m_{\alpha},\,\vert\alpha\vert n - 1 > 0,$ hence $\frac{1}{\vert\alpha n - 1\vert}\leq\frac{1}{\vert\alpha\vert n - 1}$ and $\frac{1}{\vert\alpha (n + 1)- 1\vert}\leq\frac{1}{\vert\alpha\vert n - 1}.$ Consequently $\vert\frac{n}{\alpha n - 1}-\frac{n+1}{\alpha (n+1)- 1}\vert=\frac{1}{\vert\alpha n - 1\vert\vert\alpha(n+1) - 1\vert}\leq\frac{1}{(\vert\alpha\vert n - 1)^{2}}.$ Then $\Sigma_{n=m_{\alpha}}^{\infty}\frac{n}{\alpha n - 1}e_{n}\in A$ by Proposition II.7 and therefore $\Sigma_{n\in I_{\alpha}}\frac{n}{\alpha n - 1}e_{n}\in A.$ If $\alpha^{-1}\notin \mathbb{N}\,$ i.e. $I_{\alpha}= \mathbb{N},\,(\alpha e - x)\Sigma_{n=1}^{\infty}\frac{n}{\alpha n - 1}e_{n} = e,$ then $\alpha e - x$ is invertible, a contradiction. If $\alpha^{-1}\in \mathbb{N},$ put $n_{\alpha}=\alpha^{-1},$ then $ (\alpha e - x)\Sigma_{n\in I_{\alpha}}\frac{n}{\alpha n - 1}e_{n} = \Sigma_{n\in I_{\alpha}}e_{n},$ hence $\Sigma_{n\in I_{\alpha}}e_{n}\in Ker(f)$ since $\alpha e - x\in Ker(f).$ If $e_{n_{\alpha}}\in Ker(f),$ then $e = e_{n_{\alpha}}+ \Sigma_{n\in I_{\alpha}}e_{n}\in Ker(f),  $ a contradiction. Finally $f(e_{n_{\alpha}})\neq 0$ and therefore $f = e_{n_{\alpha}}^{\ast}.$
 
\noindent \textbf{}

\noindent \textbf{Remark.} Proposition II.5 and Theorem II.8 give an answer to Sawon and Wronski's problem.

\noindent \textbf{}

\noindent \textbf{Note.} Let $A$ be a unital $B_{0}$-algebra with an orthogonal basis $(e_{n})_{n\geq 1}$ and let $ \delta $ be a translation invariant metric defining the topology of $ A. $ For each $ n\geq 1, $ define $ b^{(n)} = \Sigma_{k=n}^{\infty}e_{k}, $ this sequence converges to zero in $ A. $ For each $ n\geq 1, $ there exists $ p_{n}\geq 2\: $ such that $ \delta(b^{(k)},0)\leq\frac{1}{n^{2}} $ for all $ k\geq p_{n}. $ We can assume that $(p_{n})_{n\geq 1}$ is a strictly increasing sequence in $ \mathbb{N}. $ Then $ x = \Sigma_{n=1}^{\infty}b^{(p_{n})} $ converges in $ A $ since $ \delta(b^{(p_{n})},0)\leq\frac{1}{n^{2}} $ for all $ n\geq 1. $ Let $ k\geq 1,$ we have $ e_{k}^{\ast}(x)=0 $ for $ 1\leq k < p_{1} $ and $ e_{k}^{\ast}(x)=n $ for $ p_{n}\leq k < p_{n+1}.$ Therefore $(e_{k}^{\ast}(x))_{k\geq 1}$ is an increasing sequence of natural integers such that $ e_{k}^{\ast}(x)\rightarrow\infty.$ This shows that the condition in [3, Theorem 2.1] is always satisfied and so Theorem II.8 can be deduced from Theorem II.3.

\noindent \textbf{}

\noindent \textbf{}
   
\noindent \textbf{References}

\noindent \textbf{}

\noindent [1] D. J. H. Garling, On topological sequence, Proc. Cambridge Philos. Soc., 63 (1967), 997-1019.
 
\noindent [2] T. Husain, Orthogonal Schauder bases, Monographs and Textbooks in Pure and Applied
Mathematics, 143, Marcel Dekker, New York, 1991.

\noindent [3] Z. Sawon and Z. Wronski, Fr\'{e}chet algebras with orthogonal basis, Colloquium Mathematicum, 48 (1984), 103-110.
  
\noindent [4] W. Zelazko, Metric generalizations of Banach algebras, Rozprawy Matematyczne, 47, Warszawa, 1965.

\noindent \textbf{}  

\noindent \textbf{} 

\noindent \textbf{} Ecole Normale Sup\'{e}rieure

\noindent \textbf{} Avenue Oued Akreuch

\noindent \textbf{} Takaddoum, BP 5118, Rabat

\noindent \textbf{} Morocco
 
\noindent \textbf{} 

\noindent \textbf{} E-mail:  mohammed.elazhari@yahoo.fr

\end{document}